\documentclass[11pt]{article}
\usepackage[utf8x]{inputenc}
\usepackage{url}
\usepackage{graphicx}
\usepackage[tight,footnotesize]{subfigure}
\usepackage{indentfirst}
\usepackage[margin=25mm]{geometry}
\usepackage{amsmath,amsfonts,amssymb,amsthm,mathrsfs}

\newtheorem{theorem}{Theorem}[section]
\newtheorem{lemma}[theorem]{Lemma}

\title{Projection Onto A Simplex}
\author{Yunmei Chen and Xiaojing Ye 
\thanks{Department of Mathematics,
        University of Florida,
        PO Box 118105, Gainesville, FL 32611-8105.
        Phone (352) 392-0281. Fax (352) 392-8357.
	Email: {\tt yun,xye@ufl.edu},
        Web: http://www.math.ufl.edu/$\sim$yun,xye.}
}

\begin{document}
\maketitle

\begin{abstract}
This mini-paper presents a fast and simple algorithm to compute the projection onto
the canonical simplex $\triangle^n$. Utilizing the Moreau's identity, we show that the
problem is essentially a univariate minimization and the objective function is
strictly convex and continuously differentiable. Moreover, it is shown that
there are at most $n$ candidates which can be computed explicitly, 
and the minimizer is the only one that falls
into the correct interval.
\bigskip

\end{abstract}

\section{Introduction}
The computation of projection onto simplex appears in many imaging and statistics problems,
such as multiphase segmentation, diffusion tensor imaging, etc.
The problem can be described as follows: for any given vector $y\in\mathbb{R}^n$, the
projection of $y$ onto the simplex $\triangle^n$ is to solve the minimization problem
\begin{equation}\label{eqn:ps}
 x=\arg\min_{x\in\triangle^n}\|x-y\|,
\end{equation}
where $\|\cdot\|$ denotes the regular Euclidean norm and $\triangle^n$ is the canonical simplex 
defined by
\begin{equation}
 \triangle^n:=\left\{x=(x_1,\cdots,x_n)^T\in\mathbb{R}^n:0\leq x_i\leq 1, i=1,\cdots,n, 
\mbox{ and }\sum_{i=1}^nx_i=1\right\}.
\end{equation}
For example, $\triangle^2$ is the line segment between two points $(1,0)$ and $(0,1)$ in
$\mathbb{R}^2$, and $\triangle^3$ is the triangle in $\mathbb{R}^3$ with vertexes 
$(1,0,0)$, $(0,1,0)$ and
$(0,0,1)$. The solution is nontrivial and it does not yield an explicit form.
Early attempts include iterative projections onto affine subspaces of $\mathbb{R}^n$ \cite{M86},
and Lagrangian method that shows the optimal solution is $(y-t)_+$ for some $t$ 
followed by the iterative bisection
algorithm to find this $t$, etc. In addition, it is worth pointing out that
 there are a large number of literatures
projection onto the $\ell_1$ ball $\{x\in\mathbb{R}^n:\|x\|_1\leq1\}$, e.g.~\cite{Duchi08,Berg08}, 
which has has similar formulation and is in general easier to solve. 

In this report, we present a novel, faster and simpler numerical
algorithm to solve \eqref{eqn:ps} for any dimension $n\geq2$.

\section{Algorithm}
The minimization problem in \eqref{eqn:ps} can be rewritten as
\begin{equation}\label{eqn:psequiv}
 \min_{x\in\mathbb{R}^n}\pi(x)+\frac{1}{2}\|x-y\|^2,\quad
y=(y_1,\cdots,y_n)^T\in\mathbb{R}^n,
\end{equation}
where $\pi$ is the indicator function of $\triangle^n$ defined by
\begin{equation}\label{eqn:pi}
 \pi(x)=\left\{
\begin{array}{cl}
 0 & \mbox{if }x\in\triangle^n,\\
+\infty & \mbox{otherwise }. 
\end{array}
\right.
\end{equation}
As $\pi$ is a proper and convex function due to the fact that $\triangle^n$ is close and
convex, we know \eqref{eqn:psequiv}, and hence \eqref{eqn:ps}, has a unique solution.

In the literature, the solution to \eqref{eqn:psequiv} is also denoted by the Moreau's
proximity operator \cite{moreau65}
\begin{equation}\label{eqn:piprox}
 (I+\partial\pi)^{-1}(y):=\arg\min_{x\in\mathbb{R}^n}\pi(x)+\frac{1}{2}\|x-y\|^2,
\end{equation}
which satisfies the Moreau's identity \cite{CW05}
\begin{equation}\label{eqn:mid}
y=(I+\partial\pi)^{-1}(y)+(I+\partial\pi^*)^{-1}(y),
\end{equation}
where $\pi^*$ is the Fenchel transform of $\pi$ defined by
\begin{equation}\label{eqn:pistar0}
\pi^*(y):=\sup_{z\in\mathbb{R}^n}\langle z,y\rangle-\pi(z).
\end{equation}
Note that $\pi^*$ in \eqref{eqn:pistar0} has a simple form if $\pi$ is the indicator
function of $\triangle^n$ defined in \eqref{eqn:pi}, i.e.
\begin{equation}\label{eqn:pistar}
\begin{split}
 \pi^*(y):=&\sup_{z\in\mathbb{R}^n}\langle z,y\rangle-\pi(z)
=\sup_{z\in\triangle^n}\langle z,y\rangle\\
 =&\max_{z\in\triangle^n}\sum_{i=1}^nz_iy_i=\max_{1\leq i\leq n}\{y_i\}.
\end{split}
\end{equation}
So $\pi^*(y)$ is just the largest component of $y$.
Based on \eqref{eqn:mid}, it is suffice to compute the Moreau's proximity operator of $\pi^*$
\begin{equation}\label{eqn:pistarprox}
\begin{split}
 (I+\partial\pi^*)^{-1}(y)=&\arg\min_{z\in\mathbb{R}^n}\pi^*(z)+\frac{1}{2}\|z-y\|^2\\
=&\arg\min_{z\in\mathbb{R}^n}\left\{\max_{1\leq i\leq n}\{z_i\}+\frac{1}{2}\|z-y\|^2\right\}
\end{split}
\end{equation}
and then \eqref{eqn:piprox} can be obtained by
\begin{equation}\label{eqn:mideqiv}
 (I+\partial\pi)^{-1}(y)=y-(I+\partial\pi^*)^{-1}(y).
\end{equation}

To find the solution of \eqref{eqn:pistarprox}, we first sort the components of 
$y=(y_1,\cdots,y_n)^T\in\mathbb{R}^n$ in the ascending order as $y_{(1)}\leq\cdots\leq y_{(n)}$.
Then we know the minimization problem
in \eqref{eqn:pistarprox} can be written as 
\begin{equation}\label{eqn:pstarproxeqiv}
 \min_{z\in\mathbb{R}^n}\left\{\max_{1\leq i\leq n}\{z_i\}+\frac{1}{2}\|z-y\|^2\right\}
=\min_{t\in\mathbb{R}}\min_{z\in\mathbb{R}^n}\left\{t+\frac{1}{2}\|z-y\|^2:
\max_{1\leq i\leq n}\{z_i\}=t\right\}
\end{equation}
by introducing the new variable $t$. For any fixed $t$, the inner minimization problem on the right of
\eqref{eqn:pstarproxeqiv} is
\begin{equation}\label{eqn:innmin}
 \min_{z\in\mathbb{R}^n}\left\{t+\frac{1}{2}\|z-y\|^2:\max_{1\leq i\leq n}\{z_i\}=t\right\}
\end{equation}
and the minimizer $\hat{z}(t)$ (as it depends on $t$) is obviously
\begin{equation}\label{eqn:zhat}
\left(\hat{z}(t)\right)_i=\left\{
\begin{array}{cl}
 t & \mbox{if } y_i> t\\
 y_i & \mbox{if } y_i\leq t
\end{array}
\right.,\quad i=1,\cdots,n,
\end{equation}
and the minimum value of \eqref{eqn:innmin} is
\begin{equation}\label{eqn:f}
f(t):=t+\frac{1}{2}\|\hat{z}(t)-y\|^2=\left\{
\begin{array}{ll}
 t+\displaystyle\frac{1}{2}\sum_{j=1}^n(t-y_j)^2 & \mbox{if } t\leq y_{(1)}\\
 t+\displaystyle\frac{1}{2}\sum_{j=i+1}^n(t-y_{(j)})^2 & \mbox{if } y_{(i)}
\leq t\leq y_{(i+1)},\; i=1,\cdots, n-1\\
 t+\displaystyle\frac{1}{2}(t-y_{(n)})^2 & \mbox{if } t\geq y_{(n)}\\
\end{array}
\right.
\end{equation}
Note that $f$ is well-defined at $y_{(1)},\cdots,y_{(n)}$ as shown in Lemma \ref{lemma:f} below.
Therefore, \eqref{eqn:pstarproxeqiv} is equivalent to the univariate minimization problem
\begin{equation}\label{eqn:minf}
\min_{t\in\mathbb{R}}f(t)
\end{equation}
where $f(t)$ is defined in \eqref{eqn:f}.

Note that $f$ in \eqref{eqn:f} is piecewise quadratic, and hence is piecewise convex and smooth.
Moreover, the following lemma shows that $f\in C^1(\mathbb{R})$.

\begin{lemma}\label{lemma:f}
The function $f$ defined in \eqref{eqn:f} is continuous on $\mathbb{R}$, 
and its derivative $f'$ exists and is also
continuous on $\mathbb{R}$.
\end{lemma}
\textit{Proof}:
According to \eqref{eqn:f}, for all $i=1,\cdots,n$, there are
\begin{equation}\label{eqn:limf}
\lim_{t\to y_{(i)}^-}f(t)=y_{(i)}+\frac{1}{2}\sum_{j=i}^n(y_{(i)}-y_{(j)})^2
=y_{(i)}+\frac{1}{2}\sum_{j=i+1}^n(y_{(i)}-y_{(j)})^2=f(y_{(i)})
\end{equation}
and
\begin{equation}\label{eqn:limf}
\lim_{t\to y_{(i)}^-}f'(t)=1+\sum_{j=i}^n(y_{(i)}-y_{(j)})
=1+\sum_{j=i+1}^n(y_{(i)}-y_{(j)})=\lim_{t\to y_{(i)}^+}f'(y_{(i)}),
\end{equation}
which prove the lemma.  $\square$

Based on Lemma \ref{lemma:f}, we can obtain the derivative $f'$ as
\begin{equation}\label{eqn:fprime}
f'(t)=\left\{
\begin{array}{ll}
 1+\displaystyle\sum_{j=1}^n(t-y_j) & \mbox{if } t\leq y_{(1)}\\
 1+\displaystyle\sum_{j=i+1}^n(t-y_{(j)}) & \mbox{if } y_{(i)}\leq t\leq y_{(i+1)},\; i=1,\cdots, n-1\\
 1+\displaystyle(t-y_{(n)}) & \mbox{if } t\geq y_{(n)}\\
\end{array}
\right.
\end{equation}
Note that $f'$ is well-defined at the points $y_{(1)},\cdots,y_{(n)}$ in \eqref{eqn:fprime} 
as it is continuous. Based on the discussion above, we have the theorem that restricts the search 
in $n$ candidates and obtains the solution to \eqref{eqn:ps} using the only candidate that
falls into the correct interval.
\begin{theorem}\label{thm:ps}
 For any vector $y\in\mathbb{R}^n$, the projection of $y$ onto $\triangle^n$ as in \eqref{eqn:ps}
is obtained
by the positive part of $y-\hat{t}$:
\begin{equation}
x=(y-\hat{t})_+.
\end{equation}
where $\hat{t}$ is the only one in $\{t_i:i=0,\cdots,n-1\}$ that falls into the
corresponding interval as follows,
\begin{equation}\label{eqn:t}
 t_i:=\frac{\sum_{j=i+1}^ny_{(j)}-1}{n-i},\quad i=0,\cdots,n-1,\;
 \mbox{where }t_1\leq y_{(1)}\mbox{ and }y_{(i)}\leq t_i\leq y_{(i+1)},\; i=1,\cdots, n-1.
\end{equation}
\end{theorem}
\textit{Proof}.
Based on Lemma \ref{lemma:f}, we know $f$ is piecewise quadratic and $f\in C^1(\mathbb{R})$.
As \eqref{eqn:psequiv} has a unique minimizer, we know that 
\eqref{eqn:pistarprox} has a unique minimizer as well, and hence the optimal $t$ is unique.
Therefore, there is only one single $t$ that has vanish derivative
i.e.~$f'(t)=0$, and this $t$ is the minimizer of $\min_{t\in\mathbb{R}}f(t)$.

Now we compute the possible choices for the optimizer $\hat{t}$.
Note that $f$ is piecewisely defined, and according
to \eqref{eqn:fprime}, there are at most $n$ points that have vanish derivative. It is easy to
check that they are those shown in \eqref{eqn:t}.

However,
since the optimal $\hat{t}$ exists and is unique, we know that there is one and only one of
$\{t_i:i=0,\cdots,n-1\}$ that can be the optimal $\hat{t}$. In another words, there is only one $t_i$ that
falls into the ``correct'' interval and hence is the optimal choice of $\hat{t}$. 

Once $\hat{t}$ is obtained, we have the minimizer of \eqref{eqn:pstarproxeqiv} as $\hat{z}(\hat{t})$
where $\hat{z}(\cdot)$ is defined in \eqref{eqn:zhat}, and hence obtain the solution $x$ to 
\eqref{eqn:ps} based on the Moreau's identity \eqref{eqn:mideqiv}:
\begin{equation}\label{eqn:x}
 x_i=(y-\hat{z}(\hat{t}))_i=\left\{
\begin{array}{cl}
 y_i-\hat{t} & \mbox{if }y_i> \hat{t},\\
 0 & \mbox{otherwise }
\end{array}
\right.
\end{equation}
which implies that $x=(y-\hat{t})_+$. This completes the proof. $\square$

Based on Theorem \ref{thm:ps}, we only need to find 
the $t_i$ in \eqref{eqn:t} that falls in the corresponding interval, 
and claim it as the optimal $\hat{t}$ for \eqref{eqn:x}. This procedure is described as in Steps
2-4 in the Algorithm 1. An interpretation of such procedure is as follows. 
Suppose we start with a very large $t$,
namely $t\geq y_{(n)}$, and let $t$ go towards negative. 
First, we can see $f'(t)=1+(t-y_{(n)})\geq1>0$ if $t\geq y_{(n)}$ and hence the optimal
$t$ cannot occur in $[y_{(n)},\infty)$. As $f'$ is continuous and $f'(y_{(n)})=1>0$, we know
$f'$ is positive near $y_{(n)}^-$. As $f$ is quadratic in $[y_{(n-1)},y_{(n)})$, this also implies 
that the optimal $t$, if exists in this interval,
can only be $t_{n-1}:=y_{(n)}-1$ based on \eqref{eqn:fprime}. Due to the existence
and uniqueness of $\hat{t}$, we can surely accept $\hat{t}=t_{n-1}$ if $y_{(n-1)}\leq t_{n-1}<y_{(n)}$,
or simply $t_{n-1}\geq y_{(n-1)}$ (since obviously $t_{n-1}<y_{(n)}$).
If
$t_{n-1}<y_{(n-1)}$, we know the optimal $t$ is not achieved yet and hence $f'(y_{(n-1)})>0$ due to the
fact that $f$ is quadratic in $[y_{(n-1)},y_{(n)})$. Repeating the similar argument, we can accept
$\hat{t}=t_{n-2}:=(y_{(n-1)}+y_{(n)}-1)/2$ if $y_{(n-2)}\leq t_{n-2}<y_{(n-1)}$, or 
simply $t_{n-2}\geq y_{(n-2)}$ as we know $f(t)$ is quadratic and keeps increasing near $y_{(n-1)}$ 
in this case (hence $t_{n-2}$ cannot be equal to or larger than $y_{(n-1)}$). Based on this analysis,
we can repeat at most $n-1$ steps of such ``compute $t_i$ -- compare to $y_{(i)}$''
processes ($i=n-1,n-2,\cdots,1$) until $\hat{t}$ is found. 
If $\hat{t}$ is still not found after $n-1$ steps, 
then it must be $t_0:=(\sum_{j=1}^ny_j-1)/n$.

Completed with the last step $x=(y-\hat{t})_+$, the algorithm is summarized in Algorithm 1.
\bigskip

\noindent
\begin{table}[t]
\centering
\textbf{Algorithm 1 (\texttt{projsplx}). Projection of $y\in\mathbb{R}^n$ onto the simplex $\triangle^n$. }
\begin{enumerate}
 \item Input $y=(y_1,\cdots,y_n)^T\in\mathbb{R}^n$;
 \item Sort $y$ in the ascending order as $y_{(1)}\leq\cdots\leq y_{(n)}$, and set $i=n-1$;
 \item Compute $t_i=\frac{\sum_{j=i+1}^ny_{(j)}-1}{n-i}$. If $t_i\geq y_{(i)}$ then set 
	    $\hat{t}=t_i$ and go to Step 5, otherwise set $i\leftarrow i-1$ and redo Step 3 if $i\geq1$ 
	    or go to Step 4 if $i=0$;
 \item Set $\hat{t}=\frac{\sum_{j=1}^ny_{j}-1}{n}$;
 \item Return $x=(y-\hat{t})_+$ as the projection of $y$ onto $\triangle^n$.
\end{enumerate}
\end{table}

\section{Numerical Examples}
We manually computed the projections of several examples of $y$ in low-dimensional (2D and 3D) cases
using Algorithm 1, and the results turned out to be exact as expected. For more general cases,
it is usually nontrivial to demonstrate the correctness of Algorithm 1 numerically.

For completeness of this report, we show several numerical examples of Algorithm 1.
The algorithm is implemented in C and compiled in MATLAB (Version R2010b) using \texttt{mex} function. 
The computations
are performed on a Lenovo ThinkPad laptop with Intel Core 2 CPU at 2.53GHz, 3GB of memory, and
GNU/Linux (Kernel version 2.6.35) operating system.

\begin{figure}[t]
 \centering
\subfigure[2D case]
{\includegraphics[width=.45\textwidth]{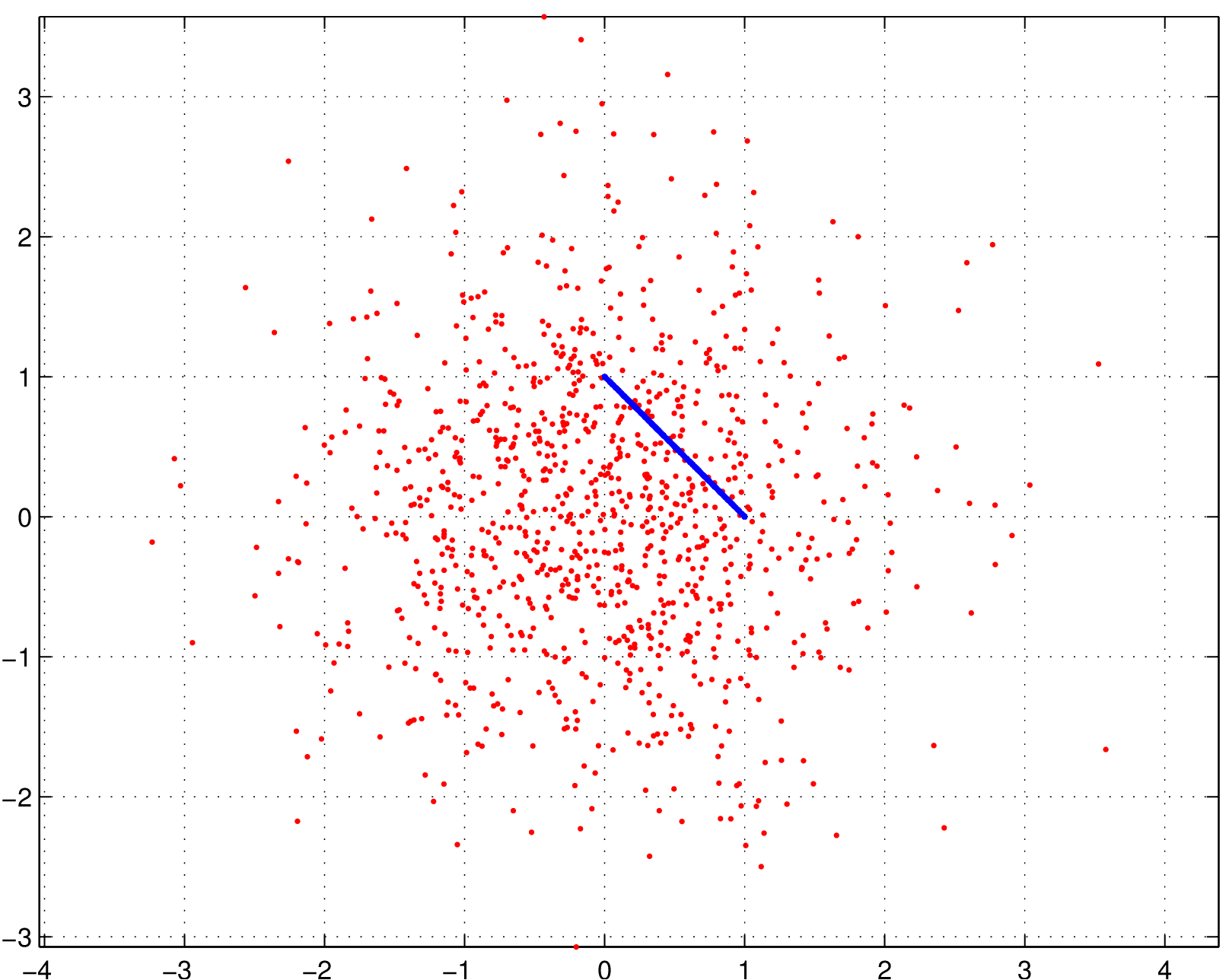}\label{fig:ps2d}}
\subfigure[3D case]
{\includegraphics[width=.45\textwidth]{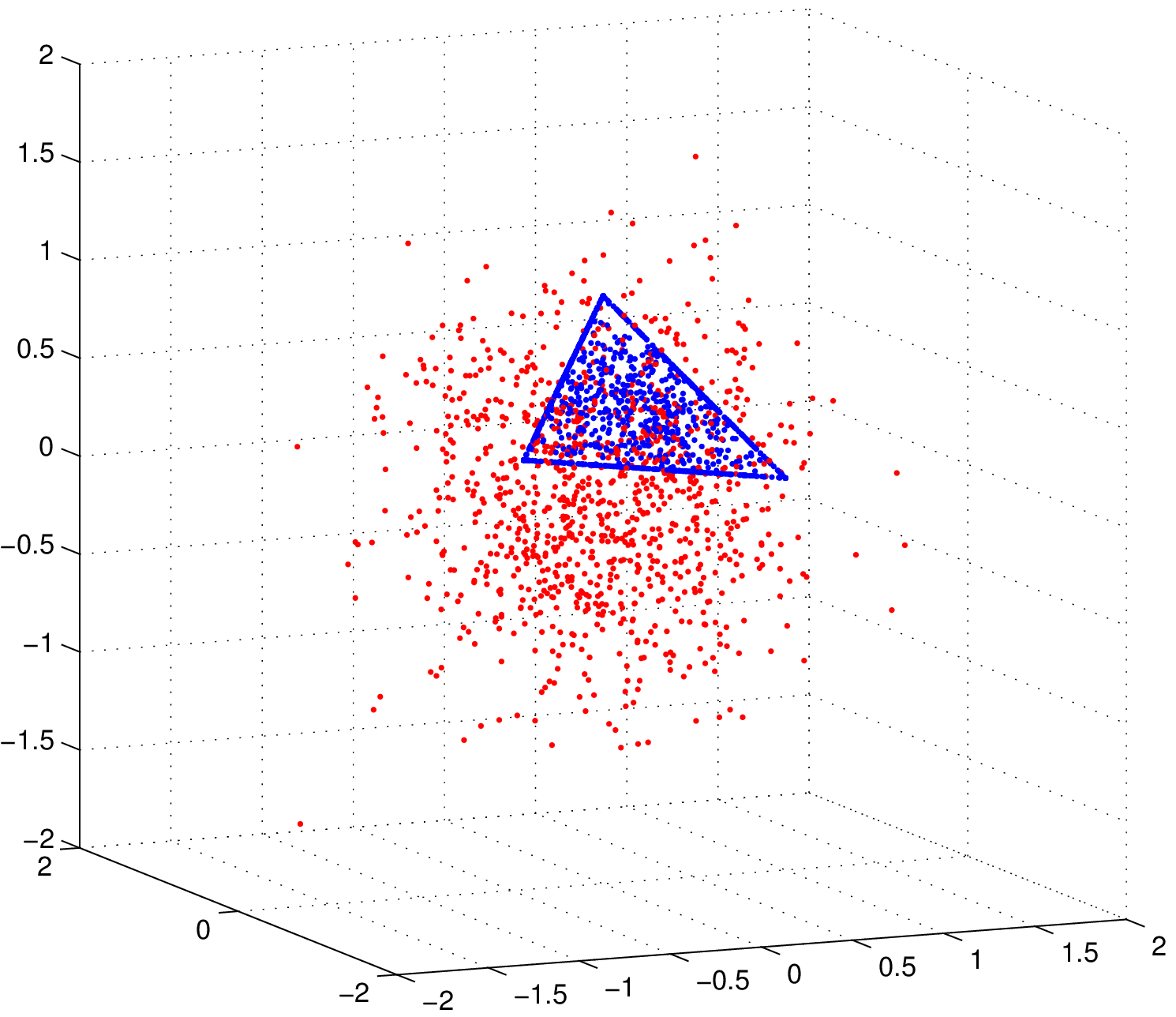}\label{fig:ps3d}}
\caption{Projections (blue) of 1024 random points (red) onto simplexes $\triangle^2$ (left) and $\triangle^3$ (right).}
\label{fig:ps}
\end{figure}

We first generate 1024 2D points from $N(0,I_2)$ using the MATLAB function \texttt{randn(2,1024)}
and project them onto $\triangle^2$ using Algorithm 1. The results
are plotted in Figure \ref{fig:ps2d}. A similar test projects 1024 samples from $N(0,0.5I_3)$ to
$\triangle^3$ and plots the results in Figure \ref{fig:ps3d}. Note that both tests are not designed to
show that Algorithm computed the correct projections, since the correspondences are not shown. 
On the other hand, these two tests demonstrate that the outputs $x$ are indeed
located at the simplex (note that we did not check anywhere in Algorithm 1 that $x$ is on the simplex).

\begin{figure}[t]
 \centering
\includegraphics[width=.8\textwidth]{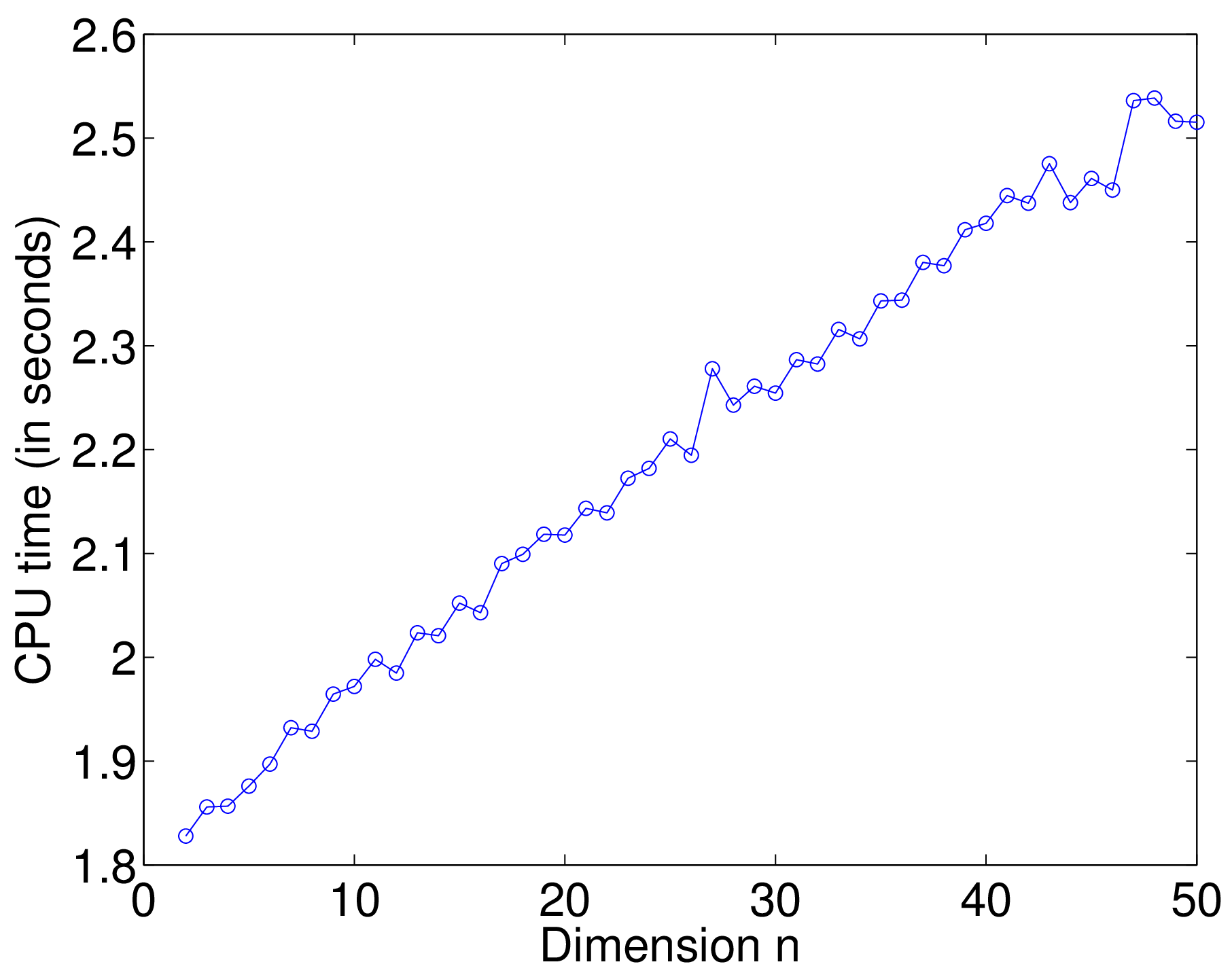}\label{fig:cpuvsn}
\caption{CPU time (in seconds) of projections of $65,536$ $n$-dimensional points drawn from
$N(0,I_n)$. The tests are carried on for $n=2,\cdots,50$.}
\label{fig:ps}
\end{figure}

The next test shows the CPU time of projections of $2^{16}=65,536$ $n$-dimensional points draw from
$N(0,I_n)$ for $n=2,3,\cdots,50$. The results are plotted in Figure \ref{fig:cpuvsn}. The CPU time
has the trend of increase as $n$ becomes larger. Interestingly, the increasing rate is rather low.
For example, for $n=5$ to $n=50$, the CPU time used only changes from 1.88s to 2.52s, for which 
the difference seems rather minor compared to the significant changes in computational complexity
of the problem.

\section{Concluding Remarks}
In this paper we propose a fast and simple algorithm \texttt{projsplx} as shown in Algorithm 1
that projects an $n$-dimensional
vector $y$ to the canonical simplex $\triangle^n$. The computation comprises of the sort of components
of the input
$y$ and at most $n$ simple ``compute-compare'' processes. The solution is exact and the algorithm
is extremely easy to implement. The MATLAB/C code is provided on the following websites for public use.

\begin{itemize}
\item Author's website: \url{http://www.math.ufl.edu/~xye/codes/projsplx.zip}
\item Matlab Central: \url{http://www.mathworks.com/matlabcentral/fileexchange/30332}
\end{itemize}

\section{Acknowledgement}
Xiaojing Ye would like to thank Stephen Becker (CalTech) for his helpful comments in complementing the references.

\bibliographystyle{plain}
\bibliography{/home/xye/Dropbox/Share/common/library}
\end{document}